\documentclass[axioms,review,accept,pdftex,moreauthors]{Definitions/mdpi}

\firstpage{1} 
\pubvolume{1}
\issuenum{1}
\articlenumber{0}
\pubyear{2024}
\copyrightyear{2024}
\externaleditor{Academic Editor: Firstname\linebreak Lastname}
\datereceived{24 May 2024} 
\daterevised{16 June 2024} 
\dateaccepted{24 June 2024} 
\datepublished{ } 
\hreflink{https://doi.org/} 

\def\Q{{\mathbb Q}}
\def\Z{{\mathbb Z}}

\def\F{{\mathbb F}}

\def\ind{{\rm{ind}}}

\Title{Monogenity and Power Integral Bases: Recent Developments}
\TitleCitation{Monogenity and Power Integral Bases: Recent Developments}

\Author{{István} 
 Gaál $^{\dagger}$\orcidA{}}
\AuthorNames{István Gaál}
\AuthorCitation{{Gaál}
, I.}
\address[1]{{Mathematical} 
 Institute, University of Debrecen, {4032 Debrecen,  Hungary}
; gaal.istvan@unideb.hu}

\firstnote{\hangafter=1 \hangindent=1.05em \hspace{-0.82em}{Current address:} 
 H--4002 Debrecen Pf.400., {Hungary.} 
}

\abstract{Monogenity is a classical area of algebraic number theory that continues to be actively researched. This paper collects the results obtained over the past few years in this area.
Several of the listed results were presented at a series of online conferences 
titled ``Monogenity and Power Integral Bases''.
We also give a collection of the most important methods used in several of these papers.
A list of open problems for further research is also given.
}



\begin{document}


\section{Introduction}

Let $K=\Q(\gamma)$ be an algebraic number field of degree $n$,
generated by the algebraic integer $\gamma$, with ring of integers $\Z_K$ and discriminant
$D_K$. It is a classical problem of algebraic number theory, going back to R. {Dedekind}
~\cite{R},
K. Hensel~\cite{hensel} and H. Hasse~\cite{hasse} to decide if the ring $\Z_K$ can be generated by
a single element $\alpha\in\Z_K$, that is, if it is mono-generated, $\Z_K=\Z[\alpha]$. In this case,
we say that the ring $\Z_K$, or the field $K$, is {\em {monogenic}}, 
 {and the integral basis} $\{1,\alpha,\ldots \alpha^{n-1}\}$
is {called a} {\em power integral basis}. 

{Recently, this area has been developing very rapidly}. In order to create a suitable forum 
to present recent results on monogenity, the author started a series of online meetings
``Monogenity and Power Integral Bases''
 (\url{https://sway.cloud.microsoft/F2kZzeZ3bmD4dFfy?ref=Link} {accessed on 15 January 2021}) 
in 2021. The purpose was 
to make contacts, circulate preprints and results, and support collaboration between
researchers all over the world working in this area. During the time of pandemy this was the only way
to contact, but later on this proved to be an easy and fast way of contacting.
Therefore, up to March 2024 we already had nine meetings and we hope to continue.

The main purpose of this paper is to give an overview of the latest developments 
in monogenity theory, about the results that were presented at the online meetings and 
the results that appeared parallel. 
The paper is also a kind of extension of the book~\cite{book} that appeared in 2019.
Most of these results are not yet contained there.

{In Section \ref{stools}, we collect the most important tools that were 
used in several works. These may be useful for further {applications}.
Section \ref{sres} collects the most important results, and finally 
in Section \ref{sfur} we try to indicate some possible
perspectives of further research.}

In favour of the reader, we collect some further concepts on monogenity.

For any primitive element $\alpha\in\Z_K$ (that is $K=\Q(\alpha)$), the {\em index}
of $\alpha$ is defined as the module index
\[
I(\alpha)=(\Z_K:\Z[\alpha]).
\]
{We} 
 obviously have 
\[
D(\alpha)=I(\alpha)^2\ D_K,
\]
where $D(\alpha)$ is the discriminant of $\alpha$,
\[
D(\alpha)=\prod_{1\le i<j\le n}(\alpha^{(i)}-\alpha^{(j)})^2,
\]
$\alpha^{(i)}$ denoting the conjugates of $\alpha$ corresponding to $\gamma^{(i)}$ ($i=1,\ldots,n$)
(in the following we shall denote similarly the conjugates of any element of $K$).
Obviously, $I(\alpha)=1$, if and only if $\Z_K=\Z[\alpha]$, that is, if 
$\{1,\alpha,\ldots \alpha^{n-1}\}$ is an integral basis, or in other words
$\alpha$ generates a power integral basis in $K$.

The {\em field index}, $m(K)$, of $K$ is defined as
\[
m(K)=\gcd \{I(\alpha)\ | \alpha\in\Z_K, K=\Q(\alpha)\}.
\]
If $K$ is monogenic, there are elements of index 1, and 
the field index is also equal to 1.
The converse is not true: the field index may happen to be 1 without the field being 
monogenic.

If $\alpha,\beta$ are primitive elements in $\Z_K$ and $\alpha+\beta\in\Z$ or $\alpha-\beta\in\Z$ then 
obviously their indices are equal. Such elements are called equivalent.
It was proved by B. J. Birch and J. R. Merriman~\cite{bime} and then in an effective
form by K. Győry~\cite{gyIII} that up to equivalence there are only finitely many 
generators of power integral bases in any number field $K$.

For any integral basis $(1,\omega_2,\ldots,\omega_n)$ of $K$ set
\[
L^{(i)}(\underline{X})=X_1+\omega_2^{(i)}X_2+\ldots+\omega_n^{(i)}X_n
\]
($i=1,\ldots,n$). Then (see~\cite{book})
\[
D(L(\underline{X}))=\prod_{1\le i<j\le n}\left(L^{(i)}(\underline{X})-L^{(j)}(\underline{X})\right)
=I(X_2,\ldots,X_n)\ D_K
\]
where $I(X_2,\ldots,X_n)$ is a homogeneous polynomial of degree $n(n-1)/2$ 
with integer coefficients, with the property that
for any primitive element $\alpha=x_1+\omega_2x_2+\ldots\omega_nx_n\in\Z_K$ we have
\[
I(\alpha)=|I(x_2,\ldots,x_n)|.
\]
The polynomial $I(X_2,\ldots,X_n)$ is called the index form corresponding to the integral basis
$(1,\omega_2,\ldots,\omega_n)$. Since equivalent algebraic integers have the same index, it
is independent of $X_1$. Therefore, determining elements $\alpha\in\Z_K$ of index $m$
is equivalent to solving the {\em index form equation}
\[
I(x_2,\ldots,x_n)=m \;\;{\rm in}\;\; x_2,\ldots,x_n\in\Z.
\]

A non-zero irreducible {\em polynomial} $f(x)\in\Z[x]$ is called {\em monogenic} if
a root $\alpha$ of $f(x)$ generates a power integral basis in the field $K=\Q(\alpha)$.
Obviously, if the polynomial $f(x)$ is monogenic, then $K$ is also monogenic,
but the converse is not true. The field $K$ may happen to be monogenic without
$f(x)$ being monogenic. The {\em index} of $f(x)$ is defined as
$\ind (f)=(\Z_K:\Z[\alpha])$.

\section{Tools}
\label{stools}

\subsection{Dedekind's Criterion}

Let  $\overline{f(x)}=\prod_{i=1}^r \overline{\phi_i(x)}^{\ell_i}$ 
modulo $p$ be the factorization of $f(x)$ modulo $p$ into powers of monic irreducible coprime polynomials of $\F_p[x]$.

For completeness, we recall here a well-known  theorem of  Dedekind:  
\begin{Theorem} [Chapter I, Proposition 8.3 of~\cite{Neu}]\

\noindent{If} 
 $p$ does not divide the index $I(\alpha)=(\Z_K:\Z[\alpha])$, then 
\[
p\Z_K=\prod_{i=1}^r \mathfrak{p}_i^{\ell_i},
\]
where $\mathfrak{p}_i=p\Z_K+\phi_i(\alpha)\Z_K$ and the residue degree of $\mathfrak{p}_i$ is 
$f(\mathfrak{p}_i)=\deg (\phi_i)$.
\label{ded1}
\end{Theorem}

As indicated above, it is very important to have a tool to determine prime divisors of 
the indices of algebraic integers. Therefore, the following well-known criterion of 
Dedekind is very frequently used:

\begin{Theorem} [Dedekind's criterion~\cite{R}, see also~\cite{Co} Theorem 6.1.4,
\cite{pz} p. 295]\

\noindent{Let}
 $f(x)\in \Z[x]$ be a monic non-zero irreducible polynomial with a root $\alpha$,
let $K=\Q(\alpha)$, and let $p$ be a prime number.
Let  $\overline{f(x)}=\prod_{i=1}^r\overline{\phi_i(x)}^{\ell_i}\mod {p}$  
be the factorization of   $\overline{f(x)}$ in $\F_p[x]$, 
with monic $\phi_i\in\Z[x]$, such that their reductions
$\overline{\phi_i(x)}$ are irreducible and pairwise coprime over $\F_p$.
Set
\[
M(x)=\cfrac{f(x)-\prod_{i=1}^r{\phi_i}^{\ell_i}(x)}{p}. 
\]
Then $M(x)\in \Z[x]$ and the following statements are equivalent:
\begin{enumerate}
\item[1.]
$p$ does not divide the index $I(\alpha)=(\Z_K:\Z[\alpha])$.
\item[2.]
For every $i=1,\dots,r$, either $\ell_i=1$ or $\ell_i\ge 2$ and $\overline{\phi_i (x)}$ does not divide 
$\overline{M(x)}$ in $\F_p[x]$.
\end{enumerate}
\end{Theorem}

\subsection{The Field Index}

We also recall a simple but very important statement of Hensel:

\begin{Theorem} [K. Hensel~\cite{hensel} p. 280]\

\noindent{The} 
 prime factors of the field index are smaller than the degree of the field.
\end{Theorem}

Denote by $\nu_p(k)$ the highest power of the prime $p$ dividing the integer $k$.

\begin{Theorem} [H. T. Engstr\"om~\cite{eng}]\

\noindent{For}  
 number fields of degree $n\le 7$, 
$\nu_p(m(K))$ is explicitly determined by the factorization of $p$
into powers of prime ideals of $p\Z_K$. 
\end{Theorem}

The corresponding tables of~\cite{eng}
are too long to include here, but they present the explicit exponents.

\subsection{Newton Polygon Method}

If $p$ divides the index $I(\alpha)=(\Z_K:\Z[\alpha])$  then Dedekind's
Theorem \ref{ded1} cannot be applied. 

Using Newton polygons, an alternative method was given by Ore~\cite{O} 
to calculate $I(\alpha)=(\Z_K:\Z[\alpha])$, $D_K$ and the prime ideal factorization of 
primes in $\Z_K$. This was further developed among others by
J. Montes and E. Nart~\cite{MN},
Fadil, L.E.  J. Montes and E. Nart~\cite{EMN} and
L. El Fadil~\cite{El}. This theory was extended to so-called higher-order Newton polygons
by J. Guardia, J. Montes and E. Nart~\cite{GMN}. The method is also called Montes algorithm.

Here we only give a short introduction to some basic notions and statement of 
this very technical method, based on the explanation used in~\cite{FG8}. 
During recent years, 
a huge amount of papers have applied this method.

For any  prime $p$, let $\nu_p$ be the $p$-adic valuation of $\Q$.
Denote by $\Q_p$ its $p$-adic completion and by $\Z_p$ the ring of $p$-adic integers. 
Let  $\nu_p$ be the Gauss's extension of $\nu_p$ to $\Q_p(x)$,
$\nu_p(P)=\min(\nu_p(a_i), \, (i=0,\dots,n)$ for any polynomial $P(x)=\sum_{i=0}^na_ix^i\in\Q_p[x]$,  
and extended by $\nu_p(P/Q)=\nu_p(P)-\nu_p(Q)$ for $0\ne P,Q\in\Q_p[x]$. Let 
$\phi\in\Z_p[x]$ be a monic polynomial whose reduction is irreducible  in
$\F_p[x]$, and let $\mathbb{F}_{\phi}$ be 
the field $\F_p[x]/(\overline{\phi})$. For any
monic polynomial  $f(x)\in \Z_p[x]$, upon  the Euclidean division
by successive powers of $\phi$, we  expand $f(x)$ as follows:
\[
f(x)=\sum_{i=0}^\ell a_i(x)\phi(x)^{i}.
\] 
This is called the $\phi$-{\em expansion} of $f(x)$  
($\deg (a_i(x))<\deg (\phi), i=1,\ldots, \ell$). 
The $\phi$-{\em Newton polygon} of $f(x)$ with 
respect to $p$ is the lower boundary convex envelope of the set of 
points $\{(i,\nu_p(a_i(x))),\, a_i(x)\neq 0\}$ in the Euclidean plane, 
which we denote by $N_{\phi}{f}$.  
The $\phi$-Newton polygon of $f$ is the process of joining the edges  
$S_1,\dots,S_r$ ordered by increasing slopes, which can be expressed 
as 
\[
\mathbb{N}_{\phi}f=S_1+\dots + S_r.
\]
For every side, $S_i$, of $\mathbb{N}_{\phi}{f}$, the {\em length} of $S_i$, denoted 
$\ell(S_i)$,  is the  length of its projection to the $x$-axis.
Its {\em height}, 
denoted by $h(S_i)$, is the  length of its projection to the $y$-axis. 
Let $d(S_i)=\gcd(\ell(S_i), h(S_i))$ be the ramification degree of $S$.
The {\it principal} $\phi$-{\em Newton polygon} of $f$,
denoted $\mathbb{N}^+_{\phi}f$, is the part of the polygon $\mathbb{N}_{\phi}f$, 
which is  determined by joining all sides of negative slopes.
For every side, $S$, of $\mathbb{N}^+_{\phi}f$, with initial point 
$(s, u_s)$ and length $\ell$, and for every 
$0\le i\le \ell$, we attach   the
{\em residue coefficient} $c_i\in\mathbb{F}_{\phi}$ as follows:
\[
c_{i}=
\displaystyle\left\{
\begin{array}{ll} 
0,& \mbox{ if } (s+i,{\it u_{s+i}}) \mbox{ lies strictly above } S,\\
\left(\dfrac{a_{s+i}(x)}{p^{{\it u_{s+i}}}}\right)
\,\,
\mod (p,\phi(x)),&\mbox{ if }(s+i,{\it u_{s+i}}) \mbox{ lies on }S,
\end{array}
\right.
\]
where $(p,\phi(x))$ is the maximal ideal of $\Z_p[x]$ generated by $p$ and $\phi$. 
Let $\lambda=-h/e$ be the slope of $S$, where  $h$ and $e$ are two positive coprime integers. Then  $d=\ell/e$ is the {\em degree} of $S$.  The points  with integer coordinates lying on $S$ are exactly 
\[
\displaystyle{(s,u_s),(s+e,u_{s}-h),\cdots, (s+de,u_{s}-dh)}.
\] 
Thus, if $i$ is not a multiple of $e$, then 
$(s+i, u_{s+i})$ does not lie in $S$, and so $c_i=0$. The polynomial
\[
f_S(y)=t_dy^d+t_{d-1}y^{d-1}+\cdots+t_{1}y+t_{0}\in\mathbb{F}_{\phi}[y],
\]
is  called the {\em residual polynomial} of $f(x)$ associated to the side $S$, 
where, for every $i=0,\dots,d$,  $t_i=c_{ie}$.

Let $\mathbb{N}^+_{\phi}{f}=S_1+\dots + S_r$ be the principal $\phi$-Newton polygon 
of $f$ with respect to $p$. We say that $f$ is a $\phi$-{\em regular} polynomial with 
respect to $p$, if  $f_{S_i}(y)$ is square free in $\mathbb{F}_{\phi}[y]$ 
for every  $i=1,\dots,r$. The polynomial $f$ is said to be  $p$-{\em regular}  
if $\overline{f(x)}=\prod_{i=1}^r\overline{\phi_i(x)}^{\ell_i}$ 
for some monic polynomials $\phi_1,\dots,\phi_t$ of $\Z[x]$, such that 
$\overline{\phi_1},\dots,\overline{\phi_t}$ are irreducible coprime 
polynomials over $\mathbb{F}_p$ and  $f$ is  a $\phi_i$-regular 
polynomial with respect to $p$ for every $i=1,\dots,t$.

 \smallskip
 Let $\phi\in\Z_p[x]$ be a monic polynomial, such that  $\overline{\phi(x)}$ 
is irreducible in $\mathbb{F}_p[x]$. 
The $\phi$-{\em index} of $f(x)$ (cf.~\cite{EMN}),
denoted by $\ind_{\phi}(f)$, is  $\deg (\phi)$ times the number of points 
with natural integer coordinates that lie below or on the polygon $\mathbb{N}^+_{\phi}{f}$,
 strictly above the horizontal axis and strictly beyond the vertical axis 
(see Figure \ref{fig1}).

\begin{figure}[H] 

\begin{tikzpicture}[x=1cm,y=0.5cm]
\draw[latex-latex] (0,6) -- (0,0) -- (10,0) ;
\draw[thick] (0,0) -- (-0.5,0);
\draw[thick] (0,0) -- (0,-0.5);
\node at (0,0) [below left,blue]{\footnotesize $0$};
\draw[thick] plot coordinates{(0,5) (1,3) (5,1) (9,0)};
\draw[thick, only marks, mark=x] plot coordinates{(1,1) (1,2) (1,3) (2,1)(2,2)     (3,1)  (3,2)  (4,1)(5,1)  };
\node at (0.5,4.2) [above  ,blue]{\footnotesize $S_{1}$};
\node at (3,2.2) [above   ,blue]{\footnotesize $S_{2}$};
\node at (7,0.5) [above   ,blue]{\footnotesize $S_{3}$};
\end{tikzpicture}
\caption{\large  $N_{\phi}^+{f}$. \label{fig1}}
\end{figure}

In the example of Figure \ref{fig1}, $ind_\phi(f)=9\times \deg (\phi)$.\\
\smallskip

Now assume that $\overline{f(x)}=\prod_{i=1}^r\overline{\phi_i(x)}^{\ell_i}$ is the
factorization of $\overline{f(x)}$ in $\mathbb{F}_p[x]$
into monic polynomials $\phi_i\in\Z[x]$, which are irreducible and pairwise 
coprime in $\mathbb{F}_p[x]$ ($i=1,\ldots,r$).

For every $i=1,\dots,r$, let  $N_{\phi_i}^+(f)=S_{i1}+\dots+S_{ir_i}$ be the principal  
$\phi_i$-Newton polygon of $f$ with respect to $p$. 
For every $j=1,\dots, r_i$,  let 
\[
f_{S_{ij}}(y)=\prod_{k=1}^{s_{ij}}\psi_{ijk}^{a_{ijk}}(y)
\]
be the factorization of $f_{S_{ij}}(y)$ in $\mathbb{F}_{\phi_i}[y]$. 
Then we have the following index theorem of Ore.

\begin{Theorem}[Theorem of Ore, see Theorem 1.7 and Theorem 1.9 in~\cite{EMN}, Theorem 3.9 in~\cite{El},\ pp. 323--325 in~\cite{MN,O}]\

 \begin{enumerate}
 \item
 We have
\[
\nu_p(\ind (f))\ge \sum_{i=1}^r \ind_{\phi_i}(f).
\]
The equality holds if $f(x)$ is $p$-regular. 
\item
If  $f(x)$ is $p$-regular, then
\[
p\Z_K=\prod_{i=1}^r\prod_{j=1}^{r_i}
\prod_{k=1}^{s_{ij}}\mathfrak{p}^{e_{ij}}_{ijk},
\]
is the factorization of $p\Z_K$ into powers of prime ideals of $\Z_K$ 
lying above $p$, where $e_{ij}=\ell_{ij}/d_{ij}$, $\ell_{ij}$ is the length of 
$S_{ij}$,  $d_{ij}$ is the ramification degree  of   $S_{ij}$, and 
$f_{ijk}=\deg (\phi_i)\times \deg (\psi_{ijk})$ is the residue degree of the prime ideal  $\mathfrak{p}_{ijk}$ over $p$.
 \end{enumerate}
 \label{ore} 
\end{Theorem}

\subsection{Algorithmic Methods}
\label{aamm}

Several of known efficient methods for the resolutions of Diophantine equations are
related to Thue equations, cf.~\cite{book}. These methods are implemented, e.g., in Magma~\cite{magma}.
Therefore, the most efficient methods for solving index form equations also reduce the 
index form equation to Thue equations.

In cubic fields, the index form equation is a cubic Thue equation, see~\cite{book}.

The below method of I. Ga\'al, A. Peth\H o and M. Pohst~\cite{gppsys,gppsim} 
reduces the index form equations in quartic fields to a cubic and some corresponding quartic Thue
equations. This method is quite often used even nowadays, and therefore we briefly present it.

Let $K= \Q( \xi )$ be a quartic number field and
$f(x)=x^4+a_1x^3+a_2x^2+a_3x+a_4 \in \Z [x]$ the minimal polynomial of $\xi$.
We represent any $\alpha\in\Z_K$ in the form
\begin{equation}
\alpha \; = \; \frac{a_{\alpha}+x\xi+y\xi^2+z\xi^3}{d},
\label{alfa4}
\end{equation}
with $a_{\alpha},x,y,z\in\Z$, and with a common denominator $d\in\Z$.
Consider the solutions of the equation
\begin{equation}
I(\alpha)= m \;\;\;\;\; (\alpha \in \Z_K )
\label{index4}
\end{equation}
for $0<m\in\Z$. We have 

\begin{Theorem} [\cite{gppsys}]\

\noindent Let $i_m=d^6m/n$, where $n=I(\xi)$.
The element $\alpha$ of (\ref{alfa4})
is a solution of (\ref{index4}),
if and only if there is a solution $(u,v)\in\Z^2$ of the cubic equation
\begin{equation}
F(u,v)=u^3-a_2u^2v+(a_1a_3-4a_4)uv^2+(4a_2a_4-a_3^2-a_1^2a_4)v^3 = \pm i_m
\label{thue34}
\end{equation}
such that $(x,y,z)$ satisfies
\begin{eqnarray}
Q_1(x,y,z)&=&x^2 -xya_1 +y^2a_2+xz(a_1^2-2a_2)+yz(a_3-a_1a_2) \nonumber \\
                &&   +z^2(-a_1a_3+a_2^2+a_4) = u,
\nonumber \\
Q_2(x,y,z)&=&y^2-xz-a_1yz+z^2a_2 = v.
\label{MN4}
\end{eqnarray}
\label{lemma24}
\end{Theorem}\vspace{-24pt}
Equation (\ref{thue34}) is either trivial to solve (when $F$ is reducible), 
or it is a cubic Thue equation. 

For a solution $(u,v)$ of (\ref{thue34}), we set
$Q_0(x,y,z)=uQ_2(x,y,z)-vQ_1(x,y,z)$.
If $\alpha$ in (\ref{alfa4}) is a solution of (\ref{index4}), then
\begin{equation}
Q_0(x,y,z)=0.
\label{q0}
\end{equation}
If $(x_0,y_0,z_0)\in\Z^3$ is a non-trivial solution of (\ref{q0}),
with, say, $z_0\ne 0$ (such a solution can be easily found, see 
L. J. Mordell~\cite{mordell}), 
then we can parametrize the solutions $x,y,z$
in the form
\begin{equation}
x=rx_0+p,\; y=ry_0+q,\; z=rz_0,
\label{pq}
\end{equation}
with rational parameters $r,p,q$. Substituting these $x,y,z$
into (\ref{q0}), we obtain an equation of the form
\[
r(c_1p+c_2q)=c_3p^2+c_4pq+c_5q^2,
\]
with integer coefficients $c_1,\ldots,c_5$. 
Multiply the equations in (\ref{pq}) by $c_1p+c_2q$
and replace $r(c_1p+c_2q)$ by $c_3p^2+c_4pq+c_5q^2$.
Further multiply the equations in (\ref{pq}) by the square of
the common denominator of $p, q$ to obtain all integer relations
(cf.~\cite{gppsim}). We divide those by $\gcd(p,q)^2$ and obtain
\begin{equation}
kx=c_{11}p^2+c_{12}pq+c_{13}q^2,\;
ky=c_{21}p^2+c_{22}pq+c_{23}q^2,\;
kz=c_{31}p^2+c_{32}pq+c_{33}q^2,\;
\label{kpqpq}
\end{equation}
with integer $c_{ij}$ and integer parameters $p,q$.
Here, $k$ is an integer parameter with the property that 
$k$ divides the $\det(C)/d_0^2$, 
where $C$ is the 3 $\times$ 3 matrix
with entries $c_{ij}$ and $d_0$ is the gcd of its entries
(cf.~\cite{gppsim}).
Finally, substituting the $x,y,z$ in (\ref{kpqpq}) into
(\ref{MN4}) we obtain
\begin{equation}
F_1(p,q)=k^2 u,\;\; F_2(p,q)=k^2 v.
\label{F12}   
\end{equation}
According to~\cite{gppsim}, at least one of the equations in 
 Equation (\ref{F12})
is a quartic Thue equation over the original number field $K$.

\section{Results}
\label{sres}

\subsection{Pure Fields, Trinomials, Quadrinomials, etc.}

There is no doubt that the Newton polygon method has been the most powerful tool during the last couple of years.
It is frequently combined with the application of Dedekind's criterion.
{While, in 2014, S. Ahmad, T. Nakahara and M. Syed~\cite{ans14} 
investigated monogenity properties of 
pure sextic fields using their subfield structure and relative monogenity,
in 2017 T. A. Gassert~\cite{g17} already used Montes algorithm to describe monogenity
of pure fields.} Note that this is only about the monogenity of the polynomials and not the
monogenity of number fields generated by a root of the polynomial (for some corrections, see
L. El Fadil~\cite{f21r}). 

Together with Newton polygons (or instead of them), Dedekind's criterion and Engström's theorem are also
often used. The following results often deal with polynomials of similar shape.
It is important to add that, especially using Newton polygons, the whole calculation must be
performed separately, even for polynomials of similar shape.

The first results investigated monogenity in pure fields (or radical extensions)
generated by a root of an irreducible binomial
of type $x^n-m$. Assuming that $m$ is squarefree, conditions were given for the monogenity (or
non-monogenity) of such pure fields, for $n=6,8,12, \ldots$, etc. A following step was to consider general exponents like $n=2^k,2^k\cdot 3^{\ell}, \ldots$, etc., and later on $n=p^k$ with a prime $p$.
For some exponents, the more complicated case of a composite $m$ was also investigated.
Here is a list of such results, for brevity indicating only the exponents considered:
\begin{itemize}
\item Z. S. Aygin and K. D. Nguyen~\cite{an23}: $n=3$;
\item  L. El Fadil~\cite{f22a}: $n=12$;
\item  L. El Fadil~\cite{f22b}: $n=18$;
\item  L. El Fadil~\cite{f22c}: $n=20$;
\item  L. El Fadil~\cite{f20}: $n=24$;
\item  L. El Fadil~\cite{f21h}: $n=36$;
\item  Fadil, L.E.  H. Ben Yakkou and J. Didi~\cite{fyd21}: $n=42$;
\item Fadil, L.E.  H. Choulli and O. Kchit~\cite{fck23}: $n=60$;
\item L. El Fadil and M. Faris~\cite{ff23}: $n=84$;
\item H. Ben Yakkou and O. Kchit~\cite{yk22}: $n=3^k$;
\item  L. El Fadil~\cite{f21t}: $n=2\cdot 3^k$;
\item  L. El Fadil~\cite{f22d}: $n=6, 2^k\cdot 3^{\ell}$;
\item Yakkou, H.B. A. Chillali and L. El Fadil~\cite{ycf21}: $n=2^k\cdot 5^{\ell}$;
\item  L. El Fadil~\cite{f22e}: $n=3^k\cdot 7^{\ell}$;
\item L. El Fadil and A. Najim~\cite{fn22}: $n=2^k\cdot 3^{\ell}$;
\item L. El Fadil and O. Kchit~\cite{fk23}: $n=2^k\cdot 7^{\ell}$;
\item  L. El Fadil~\cite{f22f}: $n=2^k\cdot 3^{\ell}\cdot 5^t$;
\item H. Ben Yakkou and L. El Fadil~\cite{yf21}: $n=p^k$;
\item  L. El Fadil~\cite{f21}: $n=6$, $m$ composite;
\item L. El Fadil and I. Gaál~\cite{FG8}: $n=8$, $m$ composite.
\end{itemize}

The exponents $n\le 9$ with a squarefree $m$ were investigated by I. Gaál and L. Remete~\cite{gr17}, which was
extended to arbitrary $m$ by L. El Fadil and I. Gaál~\cite{fg23}, \cite{FG8}.

A typical statement from this list is the following:
\begin{Theorem}[L. El Fadil and A. Najim~\cite{fn22}]\

\noindent Let $\alpha$ be a root of the irreducible polynomial $x^{2^k\cdot 3^{\ell}}-m$
with a squarefree $m$. 
If $m\not\equiv 1 \, (\bmod \; 4)$ and $m\not\equiv \pm 1 \, (\bmod \; 9)$
then $\alpha$ generates a power integral basis in $K=\Q(\alpha)$.
If $m\equiv 1 \, (\bmod \; 4)$, or $m\equiv 1 \, (\bmod \; 9)$, 
or $k=2$ and $m\equiv -1 \, (\bmod \; 9)$, then $K$ is not monogenic.
\end{Theorem}

A next step was to consider monogenity properties of number fields generated by a root
of an irreducible trinomial of type $x^n+ax^m+b$. The field index is also often determined
by using Engström's theorem. In the following list, again we only indicate the type
of trinomials considered:

\begin{itemize}
\item L. El Fadil~\cite{f24b}: $x^4+ax+b$;
\item L. El Fadil and I. Gaál~\cite{fg22}: $x^4+ax^2+b$;
\item H. Smith~\cite{s18}: $x^4 + ax + b, x^4 + c x^3 + d$;
\item L. Jones~\cite{j24c}  showed that there exist exactly three distinct monogenic trinomials of the form $x^4+bx^2+d$ with Galois $C_4$;
\item Jakhar, A. S. Kaur and S. Kumar~\cite{jkk22c}: $ x^5+ax+b$;
\item L. El Fadil~\cite{f22g}: $x^5 +ax^2 +b$;
\item L. El Fadil~\cite{f23b}: $x^5 +ax^3 +b$;
\item L. El Fadil~\cite{f23c}: $x^6 +ax +b$;
\item A. Jakhar and S. Kumar~\cite{jk22}: $x^6+ax+b$;
\item L. El Fadil~\cite{f22h}: $x^6 +ax^3 +b$;
\item L. El Fadil and O. Kchit~\cite{fk22c}: $x^6 +ax^5 +b$;
\item A. Jakhar and S. Kaur~\cite{jk23}: $x^6 + ax^m + b$;
\item R. Ibarra, H. Lembeck, M. Ozaslan, H. Smith and K. E. Stange~\cite{iloss22}: $x^n +ax+b, x^n +cx^{n-1} +d$ 
for $n=5,6$;
\item L. El Fadil and O. Kchit~\cite{fk23c}: $x^7 +ax^3 +b$;
\item H. Ben Yakkou~\cite{y22a}: $x^7 + ax^5 + b$;
\item Jakhar, A. S. Kaur and S. Kumar~\cite{jkk23d}: $x^7+ax+b$;
\item H. Ben Yakkou~\cite{y22b}: $x^8 +ax +b$;
\item H. Ben Yakkou and B. Boudine~\cite{yb23}:  $x^8 + ax + b$;
\item Jakhar, A. S. Kaur and S. Kumar~\cite{jkk23c}: $x^8 +ax^m +b$;
 \item L. Jones~\cite{j24d} considered monogenic trinomials of type $x^8+ax^4+b$ with prescribed Galois group;
\item O. Kchit~\cite{omar23}: $x^9 + ax + b$;
\item H. Ben Yakkou and P. Tiebekabe~\cite{yt22}: $x^9+ax+b$;
\item L. El Fadil and O. Kchit~\cite{fk23d}: $x^9+ax^2+b$;
\item L. El Fadil and O. Kchit~\cite{fk24e}: $x^{12}+ax^m+b$;
\item H. Ben Yakkou~\cite{y22c}: $x^{2^r}+ax^m+b$;
\item H. Ben Yakkou and L. El Fadil~\cite{yf22}: $x^n + ax + b,\;\; n=5,6,3^k,2^k\cdot 3^{\ell},2^k\cdot 3^{\ell}+1$;
\item A. Jakhar and S. Kumar~\cite{jk24} gave explicit conditions for the non-monogenity of $x^{q^s} - ax - b$;
\item A. Jakhar~\cite{jj23}: $x^{p^s} -ax^m -b$;
\item B. Jhorar and S. K. Khanduja~\cite{jk16}: $x^n + ax + b$, showed also that $f(x)=x^n - x - 1$ is
monogenic, if and only if  $|D(f)| = n^n -(n-1)^{n-1}$ is squarefree;
\item H. Ben Yakkou~\cite{y23f}: $x^n + ax^m + b,\;\; n=p^k,s\cdot p^k, 2^k\cdot 3^{\ell}$;
\item  L. El Fadil~\cite{f23e}: $x^n + ax^m + b,\;\; n=2^k\cdot 3^{\ell}$;
\item A. Jakhar~\cite{jj22}: $x^n-ax^m-b$;
\item Jakhar, A. S. Khanduja and N. Sangwan~\cite{jks16}: $x^n+ax^m+b$;
\item Jakhar, A. S. Khanduja and N. Sangwan~\cite{jks17} gave necessary and sufficient conditions 
in terms of $a,b,m,n$ for a given prime, $p$, to divide $I(\vartheta)$, where $\vartheta$ is a root of $x^n+ax^m+b$;
\item L. Jones~\cite{j22e} considered monogenic reciprocal trinomials of type $x^{2m}+Ax^m+1$;
\item L. Jones~\cite{j21c} showed that there are infinitely many primes $p$, such that $x^6 +p x^3 +1$ is monogenic with Galois group $D_6$;
\item L. Jones~\cite{j21d} showed that $x^n+x+1$ is monogenic, if and only if its discriminant is squarefree;
\item L. Jones and T. Phillips~\cite{jp18} showed that $x^n+ax+b$ is monogenic infinitely often;
\item L. Jones and D. White~\cite{jw21} found new infinite families of monogenic trinomials of type $x^n +Ax^m+B$.
\end{itemize}
A typical statement from this list is the following:

\begin{Theorem}[L. Jones and D. White~\cite{jw21}]\

\noindent Let $n \ge 2$ be an integer, with $m \ge 1$ a proper divisor of $n$. Let $t = n/m$ and let $\kappa$
denote the squarefree kernel of $m$. Let $A$ and $B$ be positive integers with $\gcd(A,B) > 1$, 
and define 
\[
D:= \frac{t^t B^{t-1} +(1-t)^{t-1} A^t}{\gcd(A,B)^{t-1}} . 
\]
If $B$ and $D$ are squarefree, and $\gcd(A,B) = 0 \; (\bmod \kappa)$, 
then $f(x) = x^n+Ax^m+B$ is monogenic. Moreover, $D(f)$ is not squarefree if $m \ge 2$.
\end{Theorem}

{The research continued into the direction considering monogenity properties} 
of quadrinomials, quintinomials, etc.,
that is, polynomials with four, five, etc., terms and the number fields generated by
a root of these polynomials:

\begin{itemize}
\item T. A. Gassert, H. Smith and K. E. Stange~\cite{gss19}: $x^4-6x^2-kx-3$; 
\item H. Ben Yakkou~\cite{y24c}: $x^4+ax^3+bx+c$; 
\item J. Harrington and L. Jones~\cite{hj} constructed new families of quartic 
polynomials with various Galois groups, which are monogenic infinitely often;
\item A. Jakhar and R. Kalwaniya~\cite{jka23}: $x^6+ax^m+bx+c$;
\item L. Jones~\cite{j24d}: $x^8+ax^6+bx^4+ax^2+1$;
\item L. Jones~\cite{j20b} constructed infinitely many monogenic polynomials 
of degree $p$ for every odd prime $p$;
\item L. Jones~\cite{j22d}: $x^p-2ptx^{p-1}+p^2t^2x^{p-2}+1$;
\item Jakhar, A. S. Kaur and S. Kumar~\cite{jkk23a}: $x^n + ax^2 + bx + c$;
\item Jakhar, A. S. Kaur and S. Kumar~\cite{jkk23b}: $x^{p^s} - ax^n - bx^m - c$;
\item A. Jakhar~\cite{j20t}: $x^n +ax^{n-1} +bx^{n-2} +c$;
\item L. Jones~\cite{j21b} constructed infinite families of reciprocal monogenic polynomials
with prescribed Galois group;
\item L. Jones~\cite{j21e} showed that if $4\le n\ge m\ge 0$ and $\gcd(n,m) = \gcd(n,k) = 1$ 
then $x^{n-m}(x+k)^m +p$ is monogenic for infinitely many primes $p$;
\item L. Jones~\cite{j19}: $x^n+A(Bx+1)^m$;
\item L. Jones~\cite{j20a}: $x^n +t\cdot g(x)$ with $n > \deg(g)$, when $g(x)$ is monic and $\deg(g) \in\{2,3,4\}$;
\item L. Jones~\cite{j22b} constructed reciprocal monogenic quintinomials
of type $x^{2^n}+Ax^{3\cdot 2^{n-2}}+Bx^{2^{n-1}}+Ax^{2^{n-2}}+1$;
\item L. Jones~\cite{j22c} considered infinite families of monogenic quadrinomials, quintinomials and 
sextinomials.
\end{itemize}

\subsection{The Relative Case}

In addition to the absolute case (extension of $\Q$), several authors considered 
monogenity problems in the relative case (extensions of an algebraic number field), 
or even similar problems in Dedekind rings. Mostly, Dedekind's criterion is used.

\begin{itemize}
\item M. E. Charkani and A. Deajim~\cite{cd12}; 
(see also A. Deajim and L. El Fadil~\cite{df21}): $x^p-m$ over number fields;
\item M. Sahmoudi and M. E. Charkani~\cite{sc23} considered relative pure cyclic extensions;
\item A. Soullami, M. Sahmoudi and O. Boughaleb~\cite{ssb21}: $x^{3^n} +ax^{3^s} -b$ over number fields;
\item O. Boughaleb, A. Soullami and M. Sahmoudi~\cite{bss23}: $x^{p^n} + ax^{p^s} - b$ over number fields;
\item H. Smith~\cite{s21} studied 
relative radical extensions;
\item S. K. Khanduja and B. Jhorar~\cite{kj16} gave equivalent versions of Dedekind's criterion in general rings;
\item S. Arpin, S. Bozlee, L. Herr and H. Smith~\cite{abs23a},~\cite{abs23b} studied 
monogenity of number rings from a modul-theoretic perspective;
\item R. Sekigawa~\cite{s22} constructed an infinite number of cyclic relative extensions 
of prime degree that are relative monogenic.
\end{itemize}

\subsection{Composite Polynomials}

Several authors considered monogenity of composites of polynomials,
composites of binomials, etc. The authors mainly use  Dedekind's criterion.

\begin{itemize}
\item J. Harrington and L. Jones~\cite{hj20} gave conditions for the monogenity of 
$(x^m - b)^n - a$, and 
 the composition of $x^n-a$ and $x^m-b$;
\item Jakhar, A. R. Kalwaniya and P. Yadav~\cite{jky24} considered monogenity of $(x^m-b)^n - a$,
and the composition of $x^n-a$ and $x^m-b$ using a refined version of the Dedekind criterion;
\item J. Harrington and L. Jones~\cite{hj21} considered monogenity of $\Phi_{p^a}(\Phi_{2^b}(x))$, 
where $\Phi_N(x)$ is the cyclotomic polynomial of index $N$;
\item L. Jones~\cite{j21a} considered monotonically stable polynomials of type $g(f^n(x))$;
\item L. Jones~\cite{j22a} constructed infinite collections of monic Eisenstein polynomials $f(x)$, such that
$f(x^{d^n})$ are monogenic for all integers $n \ge 0$ and $d>1$; 
\item L.Jones~\cite{j23b} considered monogenity of $S_k(x^p)$, where $S_k(x)=x^3-kx^2-(k+3)x-1$
the Shanks polynomial;
\item L. Jones~\cite{j23c} considered monogenity of $f(x^p)$, where $f(x)$ is the characteristic polynomial of an
$N$th order linear recurrence;
\item J. Harrington and L. Jones~\cite{hj22} gave conditions for the monogenity of $f(x^{p^n})$,
where $f(x)=x^m+ax^{m-1}+b$;
\item S. Kaur, S. Kumar and L. Remete~\cite{kkr24} considered monogenity of $f(x^k)$, where 
$f(x)=x^d+A\cdot h(x), \deg h<d$.
\end{itemize}

Let us recall a typical statement:

\begin{Theorem}[J. Harrington and L. Jones~\cite{hj21}]\
\noindent Let $a$ and $b$ be positive integers, and let $p$ be a prime. 
Then the polynomial $\Phi_{p^a}(\Phi_{2^b}(x))$ is monogenic, 
where $\Phi_N(x)$ is the cyclotomic polynomial of index $N$.
\end{Theorem}

\subsection{Connection with primes}

L. Jones~\cite{jwws22,j23a,j24a}, and J. Harrington and L. Jones~\cite{hj23} 
detected relations of monogenity of power compositional polynomials 
with properties of primes. We present here one of these statements.

For a recurrence sequence $U_0=0,U_1=1$ and $U_n=kU_{n-1}+U_{n-2}$,  
$(U_n)$ is periodic modulo any integer. Denote by $\pi_k(m)$ its period length modulo $m$. 
The prime $p$ is called a $k$-Wall--Sun--Sun prime, if 
\[
U_{\pi_k(p)}\equiv 0\; (\bmod \; p^2).
\]
 \begin{Theorem}[L. Jones~\cite{jwws22}]\
 
\noindent Let $D=k^2+4$ if $k\equiv 1\; (\bmod\; 2)$, and $D=(k/2)^2+1$ if $k\equiv 0\; (\bmod\; 2)$. 
Suppose that $k\not\equiv 0 \; (\bmod\; 4)$ and that $D$ is squarefree.
Let h denote the class number of $\Q(\sqrt{D})$. Let $s \ge 1$ be an integer, such that, 
for every odd prime divisor $p$ of $s$, $D$ is not a square modulo $p$ and $\gcd(p, hD) = 1$.
Then
\[
x^{2s^n}-k x^{s^n}-1
\]
is monogenic for all integers $n\ge 1$, if and only if 
no prime divisor of $s$ is a k-Wall--Sun--Sun prime.
\end{Theorem}

\subsection{Number of Generators of Power Integral Bases}

Some further results considered the number of non-equivalent generators
of power integral bases:

\begin{itemize}
\item M. Kang and D. Kim~\cite{kk23} considered the number of monogenic orders in pure cubic fields;
\item J. H. Evertse~\cite{ev23} considered ``rationally monogenic'' orders of number fields;
\item S. Akhtari~\cite{a20} showed that a positive proportion of cubic number fields, when ordered by their 
discriminant, are not monogenic;
\item L. Alpöge, M. Bhargava, A. Shnidman~\cite{abs20} showed that, 
if isomorphism classes of cubic fields are ordered by absolute discriminant,
then a positive proportion are not monogenic and yet have no local obstruction to being monogenic
(that is, the index form equations represent $+1$ or $-1$ mod $p$ for all primes $p$);
\item M. Bhargava~\cite{b22} proved that an  order $O$ in a quartic number field can have
at most 2760 inequivalent generators of power integral bases (and at most 182 if $|D(O)|$ is sufficiently
large). The problem is reduced to counting the number solutions of cubic and quartic Thue 
equations, somewhat analogously like described in Section \ref{aamm}, using a refined
enumeration;
\item S. Akhtari~\cite{a22} gave another proof of Bhargava's result~\cite{b22}: 
she used the more direct approach of Section \ref{aamm} and applied
sharp bounds for the numbers of solutions of cubic and quartic Thue equations.
\end{itemize}

\subsection{Miscellaneous}

In addition to the above lists, there were several further interesting statements 
achieved for monogenity. We try to recall them here.
\begin{itemize}
\item H. H. Kim~\cite{k19} showed that the number of monogenic dihedral quartic extensions with
absolute discriminant $\le X$ is of size $O(X^{3/4}(\log X)^3)$;
\item N. Khan,  S. Katayama,  T.  Nakahara and T. Uehara~\cite{kknu16} proved that the composite of a 
totally real field with a cyclotomic field of odd conductor $\ge 3$ or even $\ge$8 
has no power integral basis;
\item N. Khan, T. Nakahara and H. Sekiguchi~\cite{kns19} proved that there are 
exactly three monogenic cyclic sextic
fields of prime-power conductor, namely $\Q(\zeta_7), \Q(\zeta_9)$ and 
the maximal real subfield of $\Q(\zeta_{13})$;
\item D. Gil-Mu\u noz and M. Tinková~\cite{mt22} considered the indices of non-monogenic simplest 
cubic polynomials;
\item L. Jones~\cite{j22f} considered infinite families of monogenic Pisot (anti-Pisot) polynomials;
\item A. Jakhar and S. K. Khanduja~\cite{jk20} gave lower bounds for the $p$-index of a polynomial;
\item M. Castillo,~\cite{c22} showed, e.g., that $\Q(\alpha_n),n\ge 1$ is monogenic, where $\alpha_0=1$
and $\alpha_{n}=\sqrt{2+\alpha_{n-1}}$ for $n\ge 1$;
\item T. Kashio and R. Sekigawa~\cite{ks22} showed that  a monogenic normal cubic field is a 
simplest cubic field for some parameter;
\item F. E. Tanoé~\cite{t17} considered monogenity of biquadratic fields using a special integer basis;
\item K. V. Kouakou and F. E. Tanoé~\cite{kt17,tk21a}, and 
F. E. Tanoé and V. Kouassi~\cite{tk21} considered monogenity of triquadratic fields;
\item Aruna C. and P. Vanchinathan~\cite{av23} showed that an infinite number of so-called
exceptional quartic fields are monogenic.
\end{itemize}

\subsection{Explicit Calculations, Algorithms}
\label{algo}

The powerful methods of Dedekind's criterion and Newton polygons often decide 
about the monogenity of number fields. However, to explicitly determine all inequivalent 
generators of power integral bases one needs to perform calculations.
These algorithms usually involve Baker-type estimates, reduction methods and enumeration
algorithms, cf.~\cite{book}. There are efficient algorithms for low degree fields and
some more complicated methods for higher degree fields. Since these procedures
usually require considerable CPU time, if the number field is of high degree, or
we need information about a large number of fields, then we turn to the so-called
``fast'' algorithms for determining ``small'' solutions.
This yields a fast method 
to determine solutions of the index form equation with absolute values, say $\le$$10^{100}$. These algorithms determine all solutions
with a high probability but do not exclude extremely large solutions (which, however,
nobody has ever met).

We collect here some recent results involving explicit determination of generators of power 
integral bases.

\begin{itemize}
\item Z. Fran\u usi\'c and B. Jadrijevi\'c~\cite{jadr} calculated generators of relative
power integral bases in  a family of quartic extensions of imaginary quadratic fields;
\item I. Gaál~\cite{g20} showed that index form equations in composites of a totally real cubic field and a complex
quadratic field can be reduced to absolute Thue equations;
\item I. Gaál~\cite{g23} showed that the index form equations in composites of a totally real field and a complex
quadratic field can be reduced to the absolute index form  equations of the totally real field;
\item I. Gaál~\cite{g21} considered generators of power integral bases in fields 
generated by monogenic trinomials of type $x^6+3x^3+3a$;
\item I. Gaál~\cite{g22} considered generators of power integral bases in fields 
generated by monogenic binomial compositions of type $(x^3-b)^2+1$;
\item I. Gaál~\cite{g24} gave an efficient method to determine all 
generators of power integral bases of pure sextic fields;
\item I. Gaál and L. Remete~\cite{gr23} considered monogenity in octic fields of type $K=\Q(\sqrt[4]{a+bi})$;
\item I. Gaál~\cite{g23b} determined  ``small'' solutions of the index form equation in $K=\Q(\sqrt[6]{m})$, for 
$-5000<m<0$, such that $x^6-m$ is monogenic (1521 fields).
Experience: $\sqrt[6]{m}$ is the only generator of power  integral bases;
\item I. Gaál~\cite{g24t} determined ``small'' solutions of index form equations in $K=\Q(\sqrt[8]{m})$, 
$-5000<m<0$, such that $x^8-m$ is monogenic (2024 fields).
Experience: $\sqrt[8]{m}$ is the only generator of power  integral bases, except for $m=-1$;
\item I. Gaál~\cite{g24c} extended~\cite{fg22} on monogenity properties of trinomials of type $x^4+ax^2+b$;
\item I. Gaál~\cite{g24b} calculated generators of power integral bases in 
families of number fields generated by a root of monogenic quartic polynomials
considered in~\cite{hj}.
\end{itemize}

In~\cite{g24c,g24b}, the method described in Section \ref{aamm} was used, in 
\cite{jadr,g24t,gr23},
its relative analogue, see~~\cite{book,gprel}.

Also, here we recall some typical statements:

\begin{Theorem}[I. Gaál~\cite{g23}]\

\noindent Let $L$ be a totally real number field, $M=\Q(\sqrt{d})$, $d<0$ squarefree, assume
$\gcd(D_L,D_M)=1$.
If $\alpha$ generates a power integral basis in $K=LM$, then $\alpha=a+\beta\pm\omega$,
where $a\in\Z$, $\beta$ generates a power integral basis in $L$ and $(1,\omega)$
is integral basis in $M$.
\end{Theorem}

\begin{Theorem}[L. El Fadil and I. Gaál~\cite{fg22}]\

\noindent  Assume $a > 1, b > 1$ and $f(x)=x^4+ax^2+b$ is irreducible and monogenic. 
 If $a$, $b$ are not of type
 \[
a=\frac{u\pm 1}{v},b=\frac{u^2-1}{4v^2}
\]
for some $u, v \in  \Z, v \ne 0, u \ne 1$, then up to equivalence the root $\alpha$ of
$f(x)$ is the only generator of power integral bases in $K = \Q(\alpha)$.   
\end{Theorem}

\section{Further Research}
\label{sfur}

The above lists of results indicate what has already been already done and what is still missing.
It would be very interesting to somehow describe monogenity properties of quartic fields and maybe
quintic fields. This would require study of  quintinomials and sextinomials. 

What general exponents of binomials and trinomials
can still be considered? Is it possible to describe in general monogenity properties
of arbitrary trinomials of degree $\le 9$? 

How can one extend the available algorithms to
be able to calculate solutions of index form equations in higher degree fields?

All these and several other questions are to be answered. 
As it is seen from the above {lists}, in addition to some new ideas, often the application of
old, forgotten methods may also help.

\vspace{6pt}

\funding{This research received no external funding.}

\dataavailability{No data were used in this article.} 

\conflictsofinterest{The author declares no conflicts of interest.}

\noindent
The author has read and agreed to the published version of the manuscript.

\acknowledgments{The author is grateful to all participants of the online meetings ``Monogenity
and power integral bases''. The talks encouraged the research on monogenity.}


\begin{adjustwidth}{-\extralength}{0cm}

\reftitle{References}

\PublishersNote{}
\end{adjustwidth}
\end{document}